\documentclass[a4paper]{article}
\usepackage[utf8]{inputenc}
\usepackage{amsmath}
\usepackage{amssymb}
\usepackage{bm}
\usepackage{cite}
\usepackage{url}
\usepackage{color}
\usepackage{graphicx}
\usepackage{todonotes}
\usepackage{hyperref}

\begin{document}

  \title{A Differentiable Solver Approach to Operator Inference}
  \author{Dirk Hartmann\footnote
         {Siemens Technology, Munich (Germany),
          \href{mailto:hartmann.dirk@siemens.com}{hartmann.dirk@siemens.com}}\footnote
         {Corresponding Author}, 
          Lukas Failer\footnote
         {Siemens Digital Industries Software, Leuven (Belgium),
          \href{mailto:lukas.failler@siemens.com}{lukas.failer@siemens.com}}}
  
  \date{April 2021}
  \maketitle

  \begin{abstract}
    Model Order Reduction is a key technology for industrial applications in the context of digital twins. Key requirements are non-intrusiveness, physics-awareness, as well as robustness and usability. Operator inference based on least-squares minimization combined with the Discrete Empirical Interpolation Method captures most of these requirements, though the required regularization limits usability. Within this contribution we reformulate the problem of operator inference as a constrained optimization problem allowing to relax on the required regularization. The result is a robust model order reduction approach for real-world industrial applications, which is validated along a dynamics complex 3D cooling process of a multi-tubular reactor using a commercial software package.\\

  \end{abstract}

\section{Introduction}\label{sec:Intro}

  Digital twins are a major trend in the field of digitalization \cite{grieves2017digital,tao2018digital,jones2020characterising,hartmann2021digital}. They extend Model-Based Systems Engineering (MBSE) concepts along the complete life cycle providing novel means for enhanced decision making, e.g. during operation. While the use of model-based technology for operations is not new, it typically requires hand-crafted real-time models and the corresponding manual efforts limit applicability. Model Order Reduction (MOR) is a key technology for digital twins \cite{hartmann020model} since it allows to seamlessly transform highly accurate models, as typically used in engineering, to models of less complexity. Trading use case specific prediction cost with prediction capabilities and accuracy it allows in many cases real-time execution during operations.
  
  MOR is a well-established field, historically build on top of intrusive solutions leveraging knowledge of the numerical discretization. Since, corresponding knowledge is typically not available from commercial simulation tools, industrial adoption has been limited. However, along with the success of data-based methods such as neural-networks in the past decade also the identification of low-dimensional dynamical systems from data has shown many successful examples. This has particularly led to the formation of the new field of Scientific or Physics-Informed Machine Learning \cite{karniadakis2021physics, willcox2021imperative}. 

  Within this contribution we aim to realize a time dependent nonlinear Model Order Reduction (MOR) which leverages physics knowledge without requiring access to the underlying solver of the Full Order Model (FOM). At the same time, it should be possible to easily integrate it into any engineering workflow, such that corresponding models could be realized by simulation engineers without deep knowledge in the underlying MOR technologies.
  
  Thereby we build on the concept of Operator Inference (OI) \cite{peherstorfer2016data} and its extension to nonlinear terms using the Discrete Empirical Interpolation Method (DEIM) \cite{benner2020operator}. The concept of OI itself is built on a least-squares optimization problem, which in industrial relevant problems requires often regularization \cite{peherstorfer2020sampling,swischuk2020learning,mcquarrie2021data}. Our novel contribution is the reformulation of the OI concepts as a constrained optimization problem, which provides more accurate and robust results. Thereby, we are inspired by the recent advancements in the field of machine learning \cite{rackauckas2020universal,chen2018neural,um2020solver}.
  
  The contribution is structured as follows: In Section \ref{sec:Example} we introduce a simplified real-world industrial example which highlights the requirements for our approach as well as will be used for validation. In Section \ref{sec:MOR}, we summarize the approach taken in \cite{benner2020operator} and highlight in detail our novel extension. In Section \ref{sec:Results}, we validate our approach via numerical experiments along the example introduced in Section \ref{sec:Example}, before we conclude in Section \ref{sec:Conclusion}.
  
%
%
\section{Example}\label{sec:Example}

\begin{figure}
    \centering
    \raisebox{4.0cm}{(a)}
    \includegraphics[height=4.4cm]{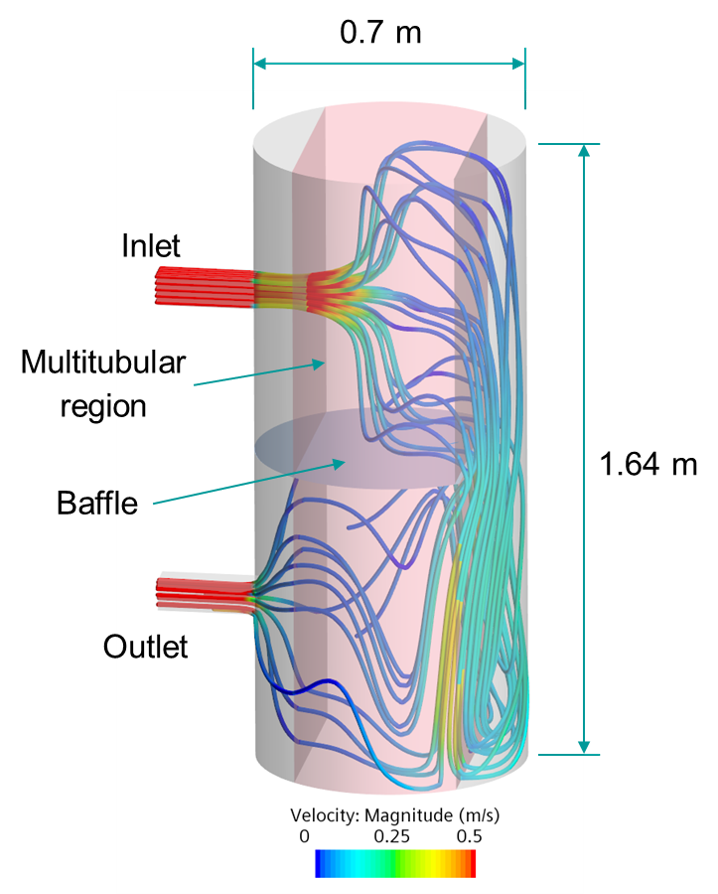}
    \raisebox{4.0cm}{(b)}
    \includegraphics[height=4.4cm]{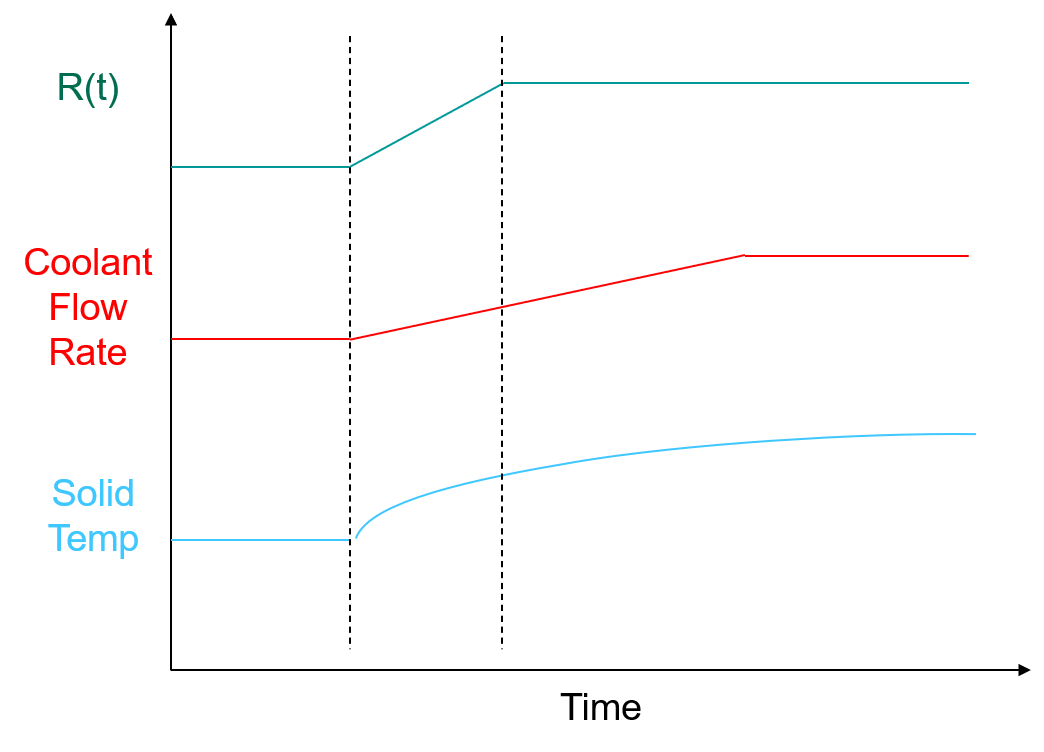}
    \caption{(a) Multi-tubular reactor where the tubular region (rose) is approximated by a porous medium model. (b) Schematic visualization of typical dynamic coolant inflow rates $v_I(t)$ (c.f. Equation \ref{eq:NS1}-\eqref{eq:NS2}), heat generation profiles $R(t)$ (c.f. Equation \eqref{eq:T04}), and resulting solid temperature at a specific spatial point.}
    \label{fig:reactor}
  \end{figure}

  Sine our contribution is motivated by fostering industrial adoption of MOR technologies, let us start with a real-world example: Goal is to realize a time dependent non-linear Reduced Order Model (ROM) for a tube and shell reactor as shown in Figure \ref{fig:reactor}. That is, the temperature distribution within the reactor shall be predicted in real-time with dynamically varying coolant inflow rates and heat generation. To achieve this, an offline-online simulation snapshot-based approach is considered. That is, snapshots of dynamic Computational Fluid Dynamics (CFD) simulations with typically many Degrees of Freedom (DoF), e.g., $> 10^5$ DoF, for a selected set of parameters are generated in an offline phase. Based on these a Reduced Order Model (ROM), e.g. $<20$ DoF, shall be developed. This model shall be capable to run in real-time to support e.g. control and operation decisions.

\subsection{Continuum Equations}\label{sec:Intro:CEq}

  The fluid velocity $\bm{v}$ and pressure field $p$ of the coolant flow within the reactor is modeled via the dynamic Navier-Stokes\footnote{In the following, we assume temperature independence. Furthermore, we would like to highlight that most commercial programs actually solve slightly more complex CFD models including sophisticated turbulence components.} equation with porous regions (c.f. Figure \ref{fig:reactor}) modeled via the Brinkman law, i.e.
  \begin{align}
    \rho_c\left(\frac{\partial\bm{v}}{\partial{t}}+(\bm{v}\cdot\nabla)\bm{v}\right) &= - \nabla p + \mu \Delta \bm{v} + \bm{\alpha}(\bm{x}) \bm{v} & &\text{in } \Omega \label{eq:NS1}\\
    \nabla\cdot\bm{v} & = 0 & &\text{in } \Omega \label{eq:NS2}
  \end{align}
  for time $t\in(0,t_\text{end})$ with temperature independent fluid density $\rho_c$, dynamic viscosity $\mu$, and porous viscous resistance $\bm{\alpha}(\bm{x})$. The flow is initiated with zero initial velocity and driven by a time-dependent inflow rate $\bm{v}_I(t)=(v_I(t),0,0)$, which is an external time-dependent parameter. Further, boundary conditions are zero flux on the reactor casing and zero normal stress at the outflow. Since the porous region $\Omega_s$ is localized, it holds $\bm{\alpha}(\bm{x})=\bm{\alpha}$ in $\Omega_s$ and $\bm{\alpha}(\bm{x})=\bm{0}$ in $\Omega \backslash \Omega_s$. 
  
  The temperature of the coolant fluid $T_c$ and the solid $T_s$ are modelled by an isotropic medium approach, i.e.
  \begin{align}
    \rho_c c_{p,c} \left(\frac{\partial T_c}{\partial{t}}+(\bm{v}\cdot\nabla)T_c\right)
       &= \nabla \cdot \left(\bm{\mathcal{K}}_c \nabla T_c\right) + \chi_{\Omega_s} {q}(T_s-T_c) & & \text{in } \Omega \label{eq:T01}\\
    \rho_s c_{p,s} \frac{\partial T_s}{\partial{t}}
       &= \nabla \cdot \left(\bm{\mathcal{K}}_s \nabla T_s\right) - {q}(T_s-T_c) + \mathcal{P}(t,T_s)& & \text{in } \Omega_s \label{eq:T02}
  \end{align}
  for time $t\in(0,t_\text{end})$ with solid density $\rho_s$, specific heat capacities $c_{p,c}$ and $c_{p,s}$, conduction coefficients $\bm{\mathcal{K}}_c(\bm{x})$ and $\bm{\mathcal{K}}_s$, volumetric heat transfer ${q}$, and heat generation $\mathcal{P}$, plus appropriate boundary and initial conditions. Here, $\chi_{\Omega_s}(\bm{x})$ represents the support of the solid region, i.e. $\chi_{\Omega_s}(\bm{x}) = 1$ for $\bm{x}\in\Omega_s$ and $\chi_{\Omega_s}(\bm{x}) = 0$ for $\bm{x}\in(\Omega\setminus\Omega_s)$. The volumetric heat transfer is is given by
  \begin{equation} \label{eq:T03}
    q(\delta T)=h\;A\;\delta T,
  \end{equation}
  where $\delta T$ is the temperature difference, $h$ is the heat transfer coefficient between the tubes and the coolant and $A$ is the exchange area density. The power is generated following an Arrhenius law
  \begin{equation} \label{eq:T04}
     \mathcal{P}(R(t),T) = {R}(t)\;\mathcal{A}\;\exp(\mathcal{B}/T),
  \end{equation}
  where ${R}(t)$ is an external dynamic parameter controlling the power output and $\mathcal{A}$, $\mathcal{B}$ are constants. 
  
  Equations \eqref{eq:NS1}-\eqref{eq:T04} constitute the Full Order Model (FOM) of the cooling process of the multi-tubular reactor. The cooling itself is influenced by the two time-dependent external parameters for the coolant inflow $v_I(t)$ and generated power $R(t)$. We assume that parameters and conditions are chosen such that the continuous FOM is well posed up to time $t_\text{end}$.

\subsection{Discretized Equations}\label{sec:Intro:DEq}

   In the following, we will mostly work with a spatially and temporally discretized version of the FOM, i.e.:
  \begin{align}
    \dot{\hat{\bm{v}}} &= \hat{\bm{A}}_1 \hat{\bm{v}}+ \hat{\bm{A}}_2 \hat{p} + \hat{\bm{H}}_1(\hat{\bm{v}} \otimes \hat{\bm{v}}) + \hat{\bm{B}}\dot{v}_I(t)  \label{eq:FOM01}\\
    0 &= \hat{\bm{A}}_3 \hat{\bm{v}}  \label{eq:FOM02}\\    
    \dot{\hat{T}}_c &= \hat{\bm{A}}_4 \hat{T}_c + \hat{\bm{A}}_5 \hat{T}_s + \hat{\bm{H}}_2(\hat{\bm{v}}\otimes\hat{T_c})  \label{eq:FOM03}\\
    \dot{\hat{T}}_s
       &= \hat{\bm{A}}_6 \hat{T}_s + \hat{\bm{A}}_7 \hat{T}_c + \hat{\bm{P}}(R(t),\hat{T}_s)  \label{eq:FOM04}
  \end{align}
  where we mark spatially discretized variables by a hat and we have introduced the term $\hat{\bm{B}}\dot{v}_I(t)$ to model time dependent inflow conditions. Taking the time derivative of the external parameterized inflow rather than the inflow itself will allow to consider the parameterization in the context of operator inference straight forward\footnote{Though within this work we will consider only q constant steady flow}. In general we cannot assume that the vectors / tensors $\hat{\bm{B}}$, $\hat{\bm{A}}_i$, $\hat{\bm{H}}_i$, and nonlinear vector function $\hat{\bm{P}}(t,\cdot)$ are known, particularly working with a commercial software package as in this contribution\footnote{3D CFD simulations are based on Simcenter STAR-CCM+ (\url{https://www.plm.automation.siemens.com/global/en/products/simcenter/STAR-CCM.html})}.
  
\section{Model Order Reduction}\label{sec:MOR}

  Within this project the discretized FOM system of time-dependent nonlinear partial differential equations (as detailed in Section \ref{sec:Intro}) shall be reduced using MOR technologies to allow for an operation parallel real-time simulation. The ROM should be significantly faster than the original model, while maintaining a sufficient degree of accuracy. At the same time, we would like it to be as interpretable and robust as possible. Thereby, the following parameters shall be able to vary: generated power controlled through $R(t)$ and coolant influx $v_I(t)$. Furthermore, the corresponding MOR technologies shall not require knowledge on the specific discretization or matrices of the discretized FOM, i.e. $\hat{\bm{B}}$, $\hat{\bm{A}}_i$, $\hat{\bm{H}}_i$, or $\hat{\bm{P}}(t,\cdot)$ in Equations \eqref{eq:FOM01}-\eqref{eq:FOM04}. 
  
  These are typical industrial requirements found in the context of MOR of large nonlinear dynamic 3D multi-physics models \cite{hartmann020model}.

  \paragraph{Approach:} Within this project, we will adopt a snapshot-based approach consisting of the two consecutive steps of reduced (latent) dimension identification and model identification in the reduced space following the approach of \cite{benner2020operator}. Thereby, we will leverage three different concepts based on optimization:
  \begin{itemize}
      \item Proper Orthogonal Decomposition (POD): Based on the set of generated snapshots, a reduced (latent) basis will be identified following \cite{chatterjee2000introduction} (reduced basis identification).
      \item Discrete Empirical Interpolation Method (DEIM): Leveraging concepts from sparse sensing we will directly reduce parts of the model, i.e. the power generation as determined by the Arrhenius law, assuming mass lumping (see e.g. \cite{thomee2007galerkin}) and following \cite{chaturantabut2010nonlinear, holmes2012turbulence, lumley1970toward} (model identification). 
      \item Operator Inference (OI): Assuming that the operators are polynomial, operator inference provides an efficient way to identify operators from given snapshot data, particularly in reduced dimensions by POD \cite{peherstorfer2016data} (model identification).
  \end{itemize}
  Following the ideas of \cite{benner2020operator} allows us to effectively address non-polynomial terms, e.g., as introduced by the Arrhenius law \eqref{eq:T04}. In the context of chemical reactors this is a crucial step, since many chemical reactions cannot be formulated in polynomial terms only, which is the basic requirement of OI adopted here \cite{qian2020lift}.
  
  Our novel contribution is to adopt a different optimization approach compared to \cite{benner2020operator}. The concept of OI requires in many cases significant regularization \cite{peherstorfer2020sampling,swischuk2020learning,mcquarrie2021data}. This clearly limits industrial applicability by non-experts. Inspired by recent works in machine learning \cite{rackauckas2020universal,chen2018neural,um2020solver} we take an additional constrained optimization approach. That is, we start with the classical least squares minimization approach adopted in OI (requiring usually regularization) and use this as an initial condition for solving another constrained optimization approach. The latter could be also interpreted as an optimal control problem, where the control parameters are the operator coefficients. Therefore we refer to this step, which is further detailed in Section \ref{sec:MOR:OC}, as operator calibration. 
  
\subsection{Proper Orthogonal Decomposition}\label{sec:MOR:POD}
  Within this approach we work with snapshot data. That is, we consider a set of $l$ different simulation time series with different control parameters $R^i(t)$ and $\dot{v}^i_I(t)$ with $1 \leq i \leq l$. Let $\hat{\bm{s}}_1^i$, $\ldots$, $\hat{\bm{s}}_k^i$ be solutions of the FOM with the corresponding control parameters $\bm{u}^i_1,\ldots\bm{u}_k^i$ at time steps $t_1$, $\ldots$, $t_k$ of the $i$th set of simulations ($1 \leq i \leq l$) where
  \begin{align} 
    \hat{\bm{s}}_j^i &= [\hat{v}_j^i, \hat{T}_{c,j}^i, \hat{T}_{s,j}^i], \label{eq:composition01} \\
         \bm{u}_j^i  &= [R_j^i,\dot{v}_{I,j}^i]. \label{eq:composition02} 
  \end{align}
  Here we have dropped the pressure $\hat{\bm{p}}_j^i$ since this is actually not a free variable but rather is a Lagrange multiplier enforcing incompressibility of the fluid.\footnote{Since we do not vary inflow conditions within our experiments in Section \ref{sec:Results}, the flow is actually a constant variable, i.e. does not change over time.}

  In our case, a single simulation snapshot $\hat{\bm{s}}_j^i$ covers more than 380 000 different values. As a first step, we therefore identify a hierarchical coordinate system, which will allow us to represent the snapshots in a much lower dimensional system. To this end we use a truncated Proper Orthogonal Decomposition (POD) or Principal Component Analysis (PCA). Let 
  \begin{equation}\label{eq:composition03}
      \bm{S}= \begin{bmatrix} 
                 \mid             &        & \mid             &        & \mid             &        & \mid             \\
                 \hat{\bm{s}}^1_1 & \ldots & \hat{\bm{s}}^1_k & \ldots & \hat{\bm{s}}^l_1 & \ldots & \hat{\bm{s}}^l_k \\
                 \mid             &        & \mid             &        & \mid             &        & \mid             \\
              \end{bmatrix}
  \end{equation}
  be the complete set of snapshot data, i.e., composed of $l$ simulation sets with different parameters and each with $k$ time steps\footnote{Since the different variables in Equation \eqref{eq:composition01} might be of different orders of magnitudes, this collection of time series data is typically re-scaled, such that all variables are of the same order of magnitude.}. That is, in total we have $m=l \cdot k$ snapshot vectors.

  The objective is now to determine a "minimal" number of modes to accurately represent the dynamics of the FOM using a Galerkin projection. Using Singular Value Decomposition (SVD) for any matrix $\bm{S}\in\mathbb{R}^{n \times m}$ (typically $m \ll n$), a truncated matrix decomposition
  \begin{equation}
      \tilde{\bm{S}} = {\bm{U}} \cdot {\bm{\Sigma}} \cdot {\bm{V}},
  \end{equation}
  of rank $r$ can be calculated. ${U} = [\psi_1,...\psi_r] \in \mathbb{R}^{n\times r}$ is a unitary matrix consisting of a truncated subset of columns of the corresponding full singular value decomposition, $\bm{\Sigma}^{r \times m}$ is a diagonal matrix with the corresponding $r$ non-negative singular values ordered from largest to smallest, and $\bm{V}^{m \times m}$ is a unitary matrix. The approach provides the best basis for $\bm{S}$ in an $l_2$ sense, i.e. $\| \bm{S}-\tilde{\bm{S}} \|$ is minimal. The specific selection of $r$ (typcially $r\ll m$ ) is an intricate manual task, but algorithmic approaches are available, e.g., using optimal thresholding \cite{gavish2014optimal}.
  
  The described method above, provides a low rank $r$-dimensional space in which we will consider the dynamics.  We would like to highlight that the method does not require explicit knowledge of the underlying dynamical system and that the matrix ${\bm{U}}$ allows to project between the high dimensional snapshot data in $\mathbb{R}^n$ and the reduced dimensional space $\mathbb{R}^r$ back and forth, i.e.
  \begin{equation}
      \hat{\bm{s}} ={\bm{U}} \cdot \bm{s} \qquad \text{and} \qquad \bm{s} ={\bm{U}}^t \cdot \hat{\bm{s}},
  \end{equation}
  where $\hat{\bm{s}}\in\mathbb{R}^n$, $\bm{s}\in\mathbb{R}^r$, and $t$ denotes the transpose. If the discretized FOM is known, the ROM can be inferred using a Galerkin projection, i.e.
  \begin{equation}\label{eq:GalerkinFOM}
      \dot{\bm{s}}(t) = {\bm{U}}^t\cdot \text{FOM}({\bm{U}}\cdot\bm{s}(t),\bm{u}(t)).
  \end{equation}
  However, in the general nonlinear case the reduction process is computationally very intensive keeping in mind that one dimension of ${\bm{U}}$ is extremely large, i.e. the evaluation of the ROM might be as challenging as the FOM. For a detailed description of dimension reduction methods and in particular POD we refer to the textbook \cite{brunton2019data}. 
  
\subsection{Discrete Empirical Interpolation Method}\label{sec:MOR:DEIM}

  Though in general we might not know the discretized FOM, for terms not involving spatial operators, e.g., integral or differential operators, we know the discretized version under the assumption of mass lumping \cite{thomee2007galerkin}, i.e., no direct coupling of degrees of freedom to other discretization points. In the case of chemical reactions, as considered here, mass lumping is often a valid assumption. This allows us to employ for some terms of the FOM \eqref{eq:FOM04} the Discrete Empirical Interpolation Method (DEIM) \cite{chaturantabut2010nonlinear,lumley1970toward,holmes2012turbulence} to identify the low dimensional equations.
  
  In our case, we want to identify a low dimensional version of the term $\hat{\bm{P}}(R(t),\bm{s})$ in \eqref{eq:FOM04}, that is to avoid evaluating 
  \begin{equation}\label{eq:GalerkinArrhenius}
    \hat{\bm{P}}(R(t),\bm{s}) = {R}(t)\;\mathcal{A}\; \tilde{\bm{U}}^T \; \exp(\mathcal{B} / ({\bm{U}}^t \cdot\bm{s})),
  \end{equation}
  and rather find a formulation not requiring $n$-times the evaluation of the exponential function plus corresponding projections. DEIM is addressing this challenge by leveraging concepts from sparse sampling and it provides an efficient evaluation of nonlinearities that scales like the rank of the POD. 
 
  DEIM is based on the snapshot matrix of the nonlinear terms:
  \begin{equation}\label{eq:DEIM01}
      \bm{P}_N = \begin{bmatrix} 
                          & \mid                                &        \\
                   \ldots & \hat{\bm{P}}(R(t),\hat{\bm{s}}^i_j) & \ldots \\
                          & \mid                                &        \\
                 \end{bmatrix},
  \end{equation}  
  where the columns are evaluations of the nonlinearity of all simulation sets $i$ at all time steps $j$. Using a corresponding low rank singular value decomposition with rank $s$\footnote{The rank $s$ can be chosen independently of $r$ but is typically of the same order.}
  \begin{equation}\label{eq:DEIM02}
      \bm{P}_N \approx \bm{U}_N \bm{\Sigma}_N \bm{V}^t_N,
  \end{equation}
  with $\bm{U}_N\in\mathbb{R}^{n \times s}$, $\bm{\Sigma}_N\in\mathbb{R}^{s \times m}$, and $\bm{V}_N\in\mathbb{R}^{m \times m}$, the method iteratively constructs a measurement matrix $\bm{P}_N \in \mathbb{R}^{n \times s}$. This measurement matrix defines the actual points where to evaluate the nonlinearity instead of evaluating it at all points in the full $n$ dimensional space of the FOM. That is, the evaluation of \eqref{eq:GalerkinArrhenius} reduces to
  \begin{equation}\label{eq:DEIM03}
    \bm{P}(R(t),\bm{s}) \approx {R}(t)\;\mathcal{A}\; \bm{P}_{1} \; \exp(\mathcal{B} / (\bm{P}_{2} \bm{s})),
  \end{equation}
  with 
  $\bm{P}_1 = {\bm{U}}^t \cdot \bm{U}_N \cdot (\bm{P}_N^t \cdot \bm{U}_N)^{-1}\in \mathbb{R}^{r \times s}$ and $\bm{P}_2 = \bm{P}^t \cdot {\bm{U}} \in \mathbb{R}^{s \times r}$. Thus, DEIM provides an efficient mean for reduction of the FOM, if the terms are known explicitly, i.e. in our case requiring the assumption of mass lumping. For a more detailed explanation of DEIM, we e.g. refer to the textbook \cite{brunton2019data}.

\subsection{(Stabilized) Operator Inference}\label{sec:MOR:OI}
  
  While we can reduce some parts of the FOM using DEIM, any parts involving spatial operators cannot be reduced by this method. To address the remaining terms, we adopt the concept of OI \cite{peherstorfer2016data, qian2020lift}. 
  
  As we have seen above, it is a fair assumption to expect that the ROM of dimension $r$ is of the same nonlinearity as the FOM of dimension $n$. Thus, in the following we assume that the ROM in the reduced space is of the form
  \begin{align}
    \dot{\bm{s}} &= \bm{A} \bm{s}+ \bm{H}(\bm{s} \otimes \bm{s}) + \bm{B}\dot{v}_I(t) + \bm{P}(R(t),\bm{s}) \label{eq:ROM01}.
  \end{align}
  That is, except the nonlinear term $\bm{P}(R(t),\cdot)$, which is given by Equation \eqref{eq:DEIM03} obtained by DEIM, it has at most quadratic polynomial terms\footnote{Since within this work we drop dependence on the coolant velocity, the model is effectively linear, i.e. the terms $\bm{H}$ and $\bm{B}$ vanish.}. 
  
  The goal of OI \cite{peherstorfer2016data, qian2020lift} is to infer the operators $\bm{A}$, $\bm{H}$, and $\bm{B}$, from a given set of time trajectory data $\hat{\bm{s}}$, respectively by its projected counterpart $\bm{s} ={\bm{U}}^t \cdot \hat{\bm{s}}$. Here, we follow the concept of \cite{qian2020lift} but introducing the additional function $\bm{P}(\cdot,\cdot) $ obtained from a DEIM following \cite{benner2020operator}. OI \cite{peherstorfer2016data} is based on the following least squares minimization problem:
  \begin{equation}\label{eq:OI_classic}
    \mathrm{argmin}_{\bm{A},\bm{H},\bm{B}} \sum_{i=1}^l\sum_{j=1}^k\big(\dot{\bm{s}}^{i}_{j} - \bm{\Theta}(\bm{s}^{i}_{j},\dot{v}_I(t),R(t))\big)^2,
  \end{equation}
  where we have introduced the following right hand side operator
  \begin{equation*}
    \bm{\Theta}(\bm{s},\bm{u},\bm{v},\bm{w}) = \bm{A}\bm{s} + \bm{H}(\bm{s}\otimes\bm{s}) + \bm{B}\dot{v}^{i}_{I,j} 
                                  + \bm{P}(R^{i}_{j},\bm{s}).
  \end{equation*}
  The optimization problem \eqref{eq:OI_classic} thereby minimizes the difference between the observed flow / time derivative of the trajectory data and the one predicted by the model \eqref{eq:ROM01} in each time step. We would like to highlight that a good estimate of the time derivatives $\dot{\bm{s}}_i$ is imperative to get optimal results (if available it is best to use exact time derivatives) \cite{peherstorfer2020sampling}.
  
  Furthermore, we would like to highlight that the Problem \eqref{eq:OI_classic} often leads to unstable dynamics, i.e. matrices might have positive eigenvalues while the underlying dynamics might not. Therefore, Problem \eqref{eq:OI_classic} typically requires appropriate regularization \cite{peherstorfer2020sampling,swischuk2020learning,mcquarrie2021data}. This can be a quite interacted task and the required manual efforts and knowledge limit industrial applicability.
  
\subsection{Operator Calibration}\label{sec:MOR:OC}

  OI poses some challenges in terms of time-derivative estimation as well as regularization. Therefore, we propose the following variant of OI: Instead of the simple quadratic least squares minimization problem \eqref{eq:OI_classic}, let us consider the corresponding constrained optimization problem: 
  \begin{align}
      &\hspace*{-2.0cm}\mathrm{argmin}_{\bm{A},\bm{H},\bm{B}} \sum_{i=1}^l\sum_{j=1}^k\left(\tilde{\bm{s}}^{i}_{j} - {\bm{s}}^{i}_{j}\right)^2 \label{eq:OI_constrained} \\
    \intertext{such that for all $1 \leq i \leq l$ \, it holds:} 
      \tilde{\bm{s}}^{i}_{j} &= \tilde{\bm{s}}^{i}_{j-1} + \delta_t  \bm{\Theta}(\tilde{\bm{s}}^{i}_{j-1},\bm{u}^{i}_{j-1}) \quad j=1,\dots,k \nonumber \\
    \tilde{\bm{s}}^{i}_{0} &=\bm{s}^{i}_{0}, \nonumber
  \end{align}
  where we have adopted an explicit Euler time discretization directly, instead of using a continuous formulation. That is, instead of optimizing the approximation of the flow / time derivatives in all time steps, the optimization problem \eqref{eq:OI_constrained} aims to minimize the difference of the complete trajectory, i.e. we take a control-like approach \cite{rawlings2017model} where the control parameters are the matrix entries. In mathematical terms, instead of solving the simple quadratic optimization \eqref{eq:OI_classic} directly with a linear solver, we have to solve a nonlinear optimization problem. Since the number of optimization parameters are rather large (number of coefficients of the matrices $\bm{A}$, $\bm{H}$, and $\bm{B}$), gradient-based optimizers are typically required. That is solving the nonlinear optimization problem \eqref{eq:OI_constrained}, we need to rely on the Karush–Kuhn–Tucker (KKT) formulation. In each optimization step, we need to solve the forward in time problem (the ROM to be identified) as well as the dual or adjoint problem corresponding to a backward in time problem. Based on the calculation of the two, the corresponding gradient for the optimization set can be determined. The corresponding scheme is independent of the specific (gradient-based) optimizer (within this work we will use MATLAB's \texttt{fminunc}) and for more details on corresponding constrained optimization we refer e.g. to \cite{rawlings2017model}. Finally, solving a nonlinear optimization problem, a good initial guess is essential. In our case we rely on the result of the classical OI as described above. 
 
  Though the problem formulation and corresponding algorithms are more complex, the novel approach provides a number of advantages. First, one can expect that less or even no regularization is required. Furthermore, a highly accurate estimation of time-derivatives is less crucial since the actual trajectories are fitted and not the time derivatives themselves. Furthermore, the approach allows to have full control of the matrices, e.g., symmetries can be enforced straight forward.
  
  We would like to highlight, that the approach is actually similar to the concepts proposed in the machine learning community \cite{chen2018neural,rackauckas2020universal,um2020solver} in the context of neural networks. These have shown a superior performance over unconstrained optimization problem. Therefore, we expect the same in our situation.
  
\section{Results}\label{sec:Results}

  As a test case we will consider a multi-tubular reactor as introduced in Section \ref{sec:Example} (c.f. Figure \ref{fig:reactor}), inspired by a real-world industrial example. Goal of the model order reduction is to provide a real-time capable model which can run in parallel to the operation and predict temperatures in the reactor, e.g., as input for a control. For simplicity, we consider here only a variation of the heat input $R(t)$, the inflow rate of the cooling fluid is kept fixed.
   
  We consider 5 different transient cases varying in the amount of heat input. For all cases we consider a coolant fluid with density $\rho_c=723 kg/m^3$, a dynamic viscosity $\mu=0.0008 Pa s$, a specific heat capacity $c_{p,c}=2590 J/(kg K)$ and a thermal conductivity $\mathcal{K}_c = 0.132 W/(mK)$. The porous solid containing the tubes has a density of  $\rho_s=3062 kg/m^3$, a specific heat capacity of $c_{p,s}=2000 J/(kg K)$, and a thermal conductivity $\mathcal{K}_c = 0.2 W/(mK)$. The porous medium has the vicious resistance $\bm{\alpha}=[8000,8000,8] kg/(m^3 s)$ and the heat exchange area density is $A=0.18337 m^2/m^3$. For all cases, we consider the following boundary and initial conditions: a coolant inflow rate of $3.0 kg/s$, initial constant temperatures in the reactor for the solid and fluid part of $533.15^\circ K$, and an inflow temperature of the coolant flow of $533.15^\circ K$. 
  
\begin{figure}
    \centering
    \includegraphics[width=12.0cm]{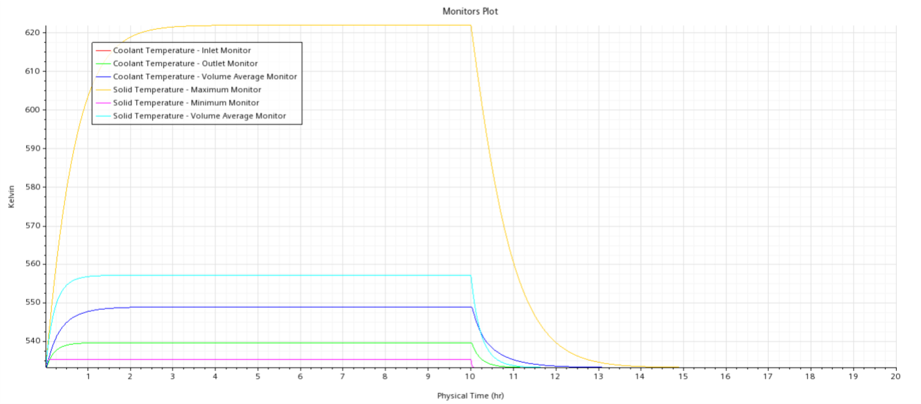}
    \caption{Heating scenario considered as a test case ($R=1$)}
    \label{fig:FOMscenario}
  \end{figure}

  The heat input is varied across the different cases as follows:  For all 5 cases,  the heating is kept constant with heat load $R$ from $0s$ till $36010s$ ($10h$) and switched off afterwards (c.f. Figure \ref{fig:FOMscenario}). Thereby, we select $R\in{0.5, 1.0, 1.5}$ as training cases, i.e. the ones used to construct the ROM, and $R\in{0.75, 1.25}$ as validation cases. We furthermore choose $\mathcal{A}=5000 J/m^3$ and $\mathcal{B} = 1500 K$.

  The corresponding FOM simulations using the geometry shown in Figure\ref{fig:reactor} are performed in the commercial CFD solver Simcenter STAR-CCM+\footnote{\url{https://www.plm.automation.siemens.com/global/en/products/simcenter/}}. The reduced order model is computed using \texttt{MATLAB}. 

\begin{figure}
    \centering
    \raisebox{3.5cm}{\small(a)}
    \includegraphics[width=5.4cm]{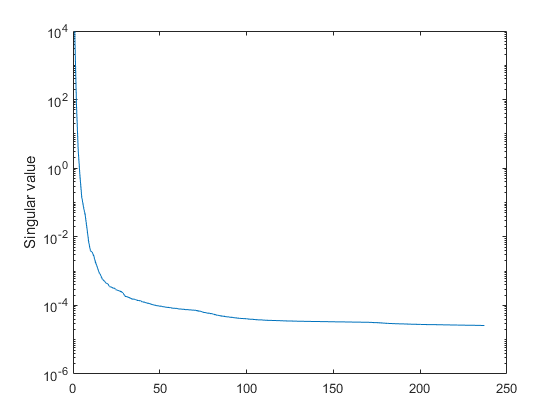}
    \raisebox{3.5cm}{\small(b)}
    \includegraphics[width=5.4cm]{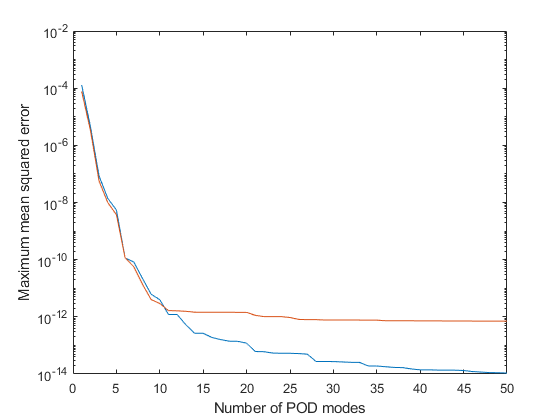}
    \caption{Proper Orthogonal Decomposition of the full 3D simulation snapshots: (a) Spectrum and (b) maximum mean squared errors (blue: training data; orange validation data).}
    \label{fig:POD_error}
  \end{figure}

  First of all let us identify the right dimension for the low dimensional model (c.f. Section \ref{sec:MOR:POD}). In Figure \ref{fig:POD_error}, the spectrum of the POD as well as the errors for the training and validation set using different number of modes are shown. For the remainder of this section we will choose 8 modes for the POD as well as for the DEIM (c.f. Section \ref{sec:MOR:DEIM}). 

\begin{figure}
    \centering
    \raisebox{3.5cm}{\small(a)}
    \includegraphics[width=5.4cm]{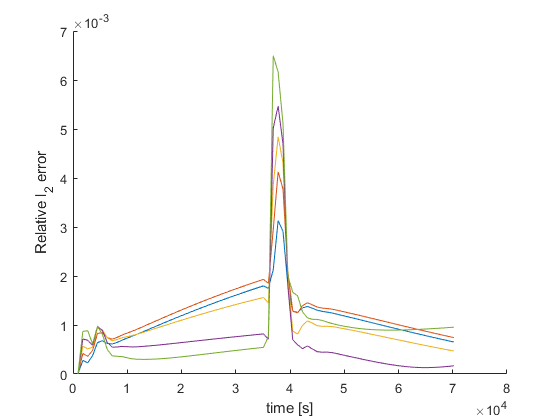}
    \raisebox{3.5cm}{\small(b)}
    \includegraphics[width=5.4cm]{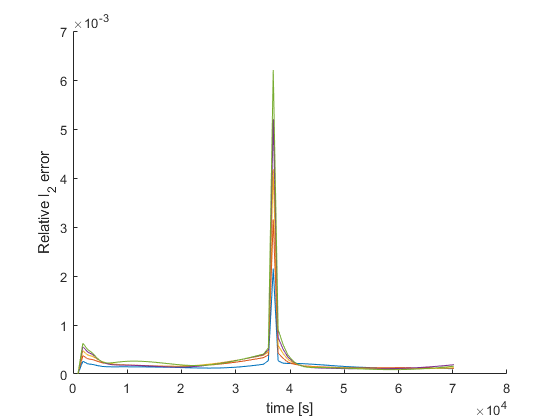} \\
    \caption{Operator Inference plus DEIM using 8 modes each: (a) Relative mean squared error of the dynamics predicted using stabilized operator inference (with stabilization parameter $\lambda= 1.0$) and (b) the same error after additional operator calibration (all 5 data sets, encoded in different color).}
    \label{fig:OI_error}
  \end{figure}

  We project the dynamics to the low dimensional space and identify in this low dimensional space the corresponding dynamical system using a combination of DEIM, OI, as well as operator calibration as described in Section \ref{sec:MOR}. This results in the following reduced order model as defined by Equations \eqref{eq:ROM01} and \eqref{eq:DEIM03}:
  \begin{equation}\label{eq:explicitROM01}
     \dot{\bm{s}} = \bm{A} \bm{s} + {R}(t)\;\mathcal{A}\; \bm{P}_{1} \; \exp(\mathcal{B} / (\bm{P}_{2} \bm{s})),
  \end{equation} 
  with the $\mathbb{R}^{8\times8}$ matrices $\bm{A}$, $\bm{P}_1$, and $\bm{P}_2$ as provided in Appendix \ref{app:ROMmatrices} (considering a POD and DEIM with 8 modes) and given constants $\mathcal{A}$ and $\mathcal{B}$ (c.f. above). Please note that we have dropped the flow of the coolant as it is constant and thus model \eqref{eq:explicitROM01} has only a linear term in addition to the Arrhenius law as mentioned above. The corresponding errors comparing the full 3D simulations projected onto the reduced low dimensional space and the simulated low dimensional ROM \eqref{eq:explicitROM01} obtained through MOR are shown in Figure \ref{fig:OI_error}. We clearly see that using OI (Tikhonov regularization with $\lambda=1.0$, for lower values the trajectories become unstable) plus operator calibration reduced the errors by one order of magnitude, except at the point in time where we switch off the heating. 
  
\begin{figure}
    \centering
    \raisebox{3.4cm}{\small(1-S)}
    \includegraphics[width=5.0cm]{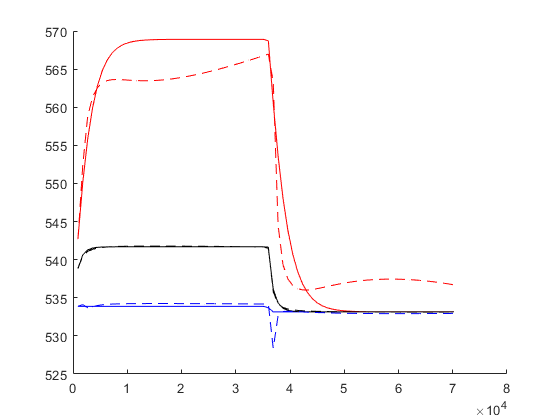}
    \raisebox{3.4cm}{\small(1-C)}
    \includegraphics[width=5.0cm]{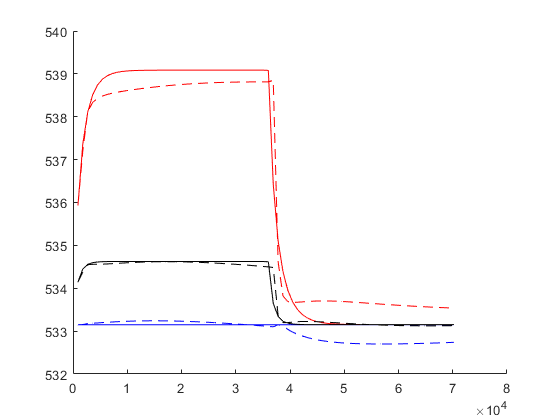}\\
    \raisebox{3.4cm}{\small(2-S)}
    \includegraphics[width=5.0cm]{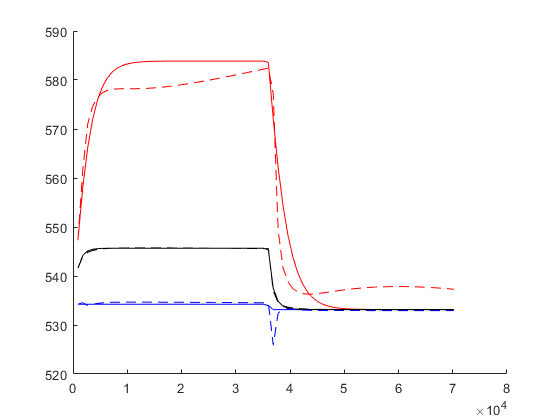}
    \raisebox{3.4cm}{\small(2-C)}
    \includegraphics[width=5.0cm]{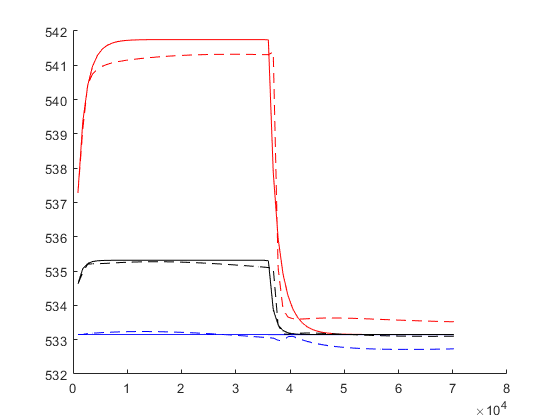}\\
    \raisebox{3.4cm}{\small(3-S)}
    \includegraphics[width=5.0cm]{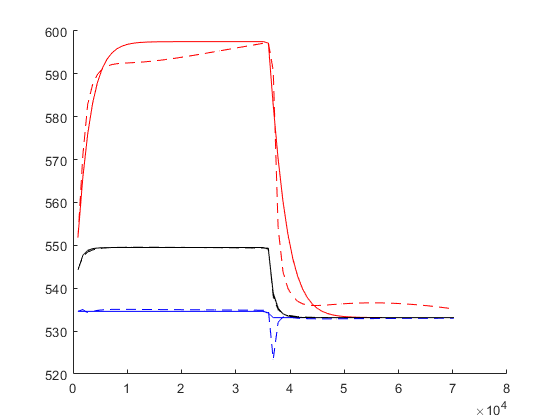}
    \raisebox{3.4cm}{\small(3-C)}
    \includegraphics[width=5.0cm]{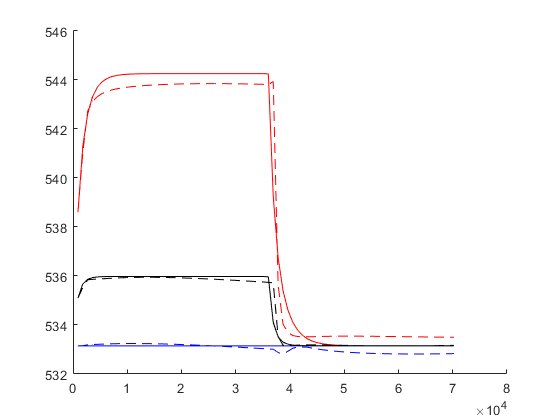}\\
    \raisebox{3.4cm}{\small(4-S)}
    \includegraphics[width=5.0cm]{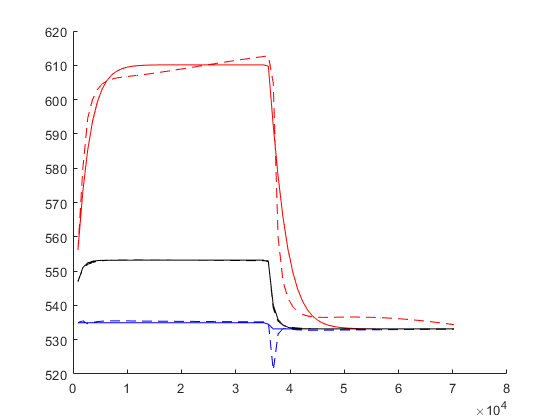}
    \raisebox{3.4cm}{\small(4-C)}
    \includegraphics[width=5.0cm]{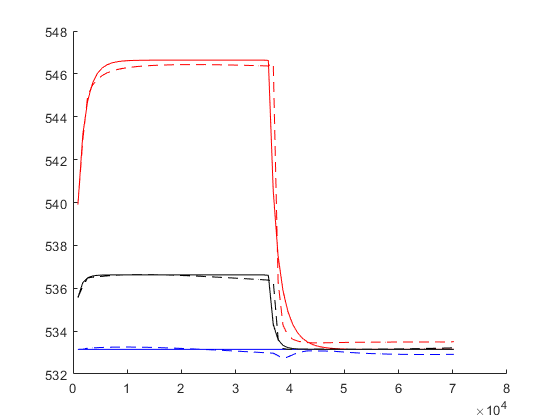}\\
    \raisebox{3.4cm}{\small(5-S)}
    \includegraphics[width=5.0cm]{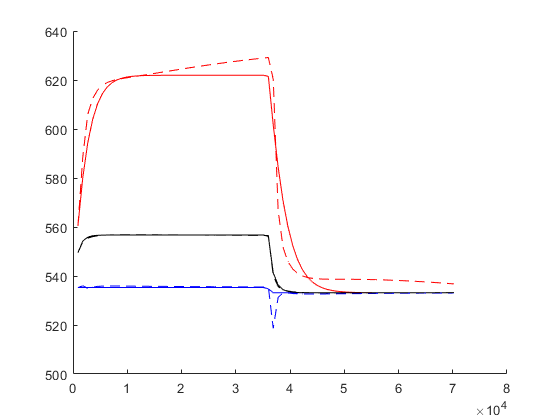}
    \raisebox{3.4cm}{\small(5-C)}
    \includegraphics[width=5.0cm]{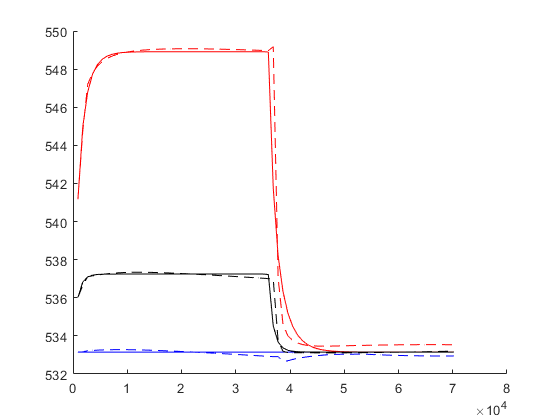}\\
    \caption{Operator Inference plus DEIM: Minimum, average, and maximum, temperatures $T_s$ ($i$-S) and $T_c$ ($i$-C) for all 5 sets ($i=1\ldots5$) over time. Reduced data is shown in dashed and full data in solid lines (using 8 POD and DEIM modes).}
    \label{fig:ROMprediction_curves}
  \end{figure}

\begin{figure}
    \centering
    \raisebox{4.6cm}{(a)}
    \includegraphics[width=10.0cm]{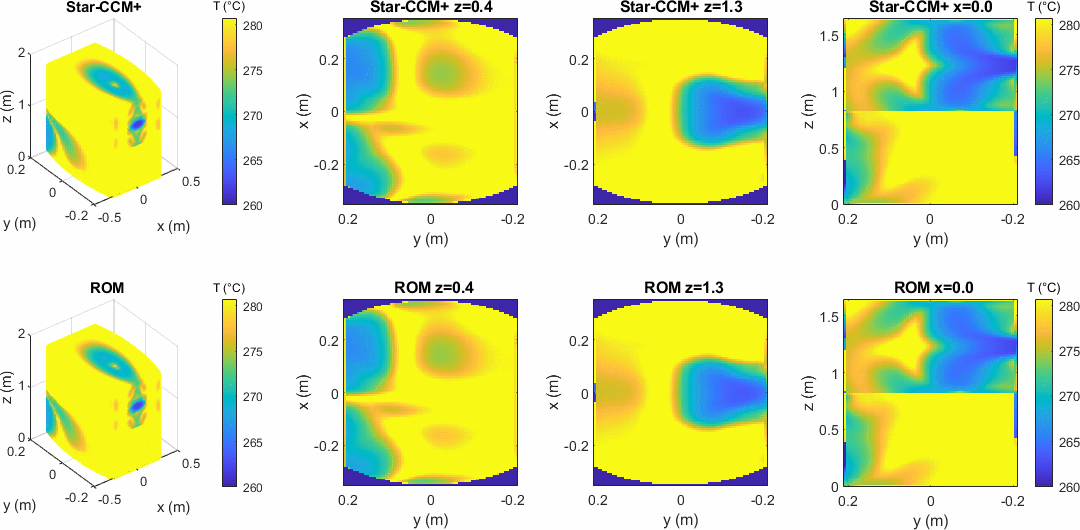}\\[0.4cm]
    \raisebox{4.6cm}{(b)}
    \includegraphics[width=10.0cm]{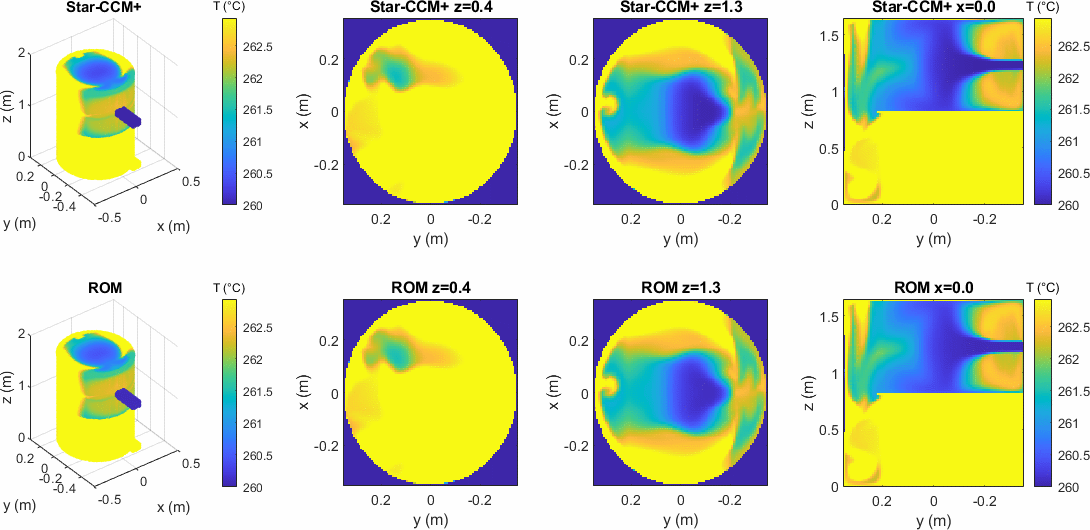}
    \caption{Comparison of spatial temperature profiles for the validation case $R=0.75$: (a) solid temperatures and (b) fluid temperatures for different planes cutting the reactor. The top row shows the original STAR-CCM+ results and the bottom row the MOR-predictions projected back into the full space.}
    \label{fig:ROMprediction_3d}
  \end{figure}

  Finally, we project the reduced order model predictions back to the full space and compare the results with the original data. On the one hand, we compare the time evolution of the average, minimum, and maximum temperature profiles of the cooling fluid as well as solid temperatures for the five different cases in Figure \ref{fig:ROMprediction_curves}. On the other hand, we also show the corresponding 2D temperature results (different cuts with planes) of the full 3D simulation for the first validation case ($R=0.75$) in Figure \ref{fig:ROMprediction_3d}. We can observe a good match of the spatial profiles as well as the average temperatures. Maximum and minimum temperatures could be further improved but are sufficiently accurate for most industrial applications.

  Performing a separate POD for the solid and fluid part might improve the predictions as suggested by \cite{benner2020operator}. In particular since the solid temperature cover a much smaller domain (only the rose part in Figure \ref{fig:reactor}). This will be further investigated in the future.

\section{Conclusion}\label{sec:Conclusion}

  Within this contribution we have addressed a combination of Proper Orthogonal Decomposition (POD), Discrete Empirical Interpolation Methods (DEIM), and Operator inference (OI) to infer a fast low dimensional model with real-time capability for operations parallel simulation following \cite{benner2020operator}. Thereby we have extended the approach introducing the concept of operator calibration, inspired by recent works in the field of machine learning \cite{rackauckas2020universal,chen2018neural,um2020solver}. 
  
  The proposed method has been validated along the case of a tubular reactor, inspired by a real-world industrial application. In particular we would like to highlight that the reduced model \eqref{eq:explicitROM01} itself is extremely compact and well suited for integration into other applications e.g. using the Functional Mock-up Interface (FMI) standard \cite{blockwitz2012functional} for model exchange or model integration in other tools.
  
  The approach is well suited for industrial applications since it allows to derive reduced order models for complex 3D multi-physics simulations without knowing any solver details. A key requirement working with commercial simulation tools. Furthermore, through the use of the novel concept of operator calibration the approach relaxes the requirements of exact estimation of time derivatives (also often not available) as well as the choice of regularization methods as required by state of the art approaches for real-world problems \cite{peherstorfer2020sampling,swischuk2020learning,mcquarrie2021data}.
  
  In future work, we plan the to extend the number of scenarios considered, e.g., including varying cooling flow and more dynamic heating profiles. Furthermore, we plan to address more complex examples as well as to explore the use of separated reduced bases for the flow and the different temperature distributions as suggested by \cite{benner2020operator}. Last but not least, the quantification of uncertainties is crucial for industrial applications and we are currently investigating corresponding concepts along the work of \cite{soize2017nonparametric,farhat2018stochastic}.

  Despite these future ambitions, we believe that the current results already show the potential of this extended concept of OI. The combined operator inference and calibration approach will likely become a standard workhorse for industrial applications in the context of digital twins.

  \paragraph*{Acknowledgements:} We would like to acknowledge the funding through the Siemens Digital Industries Software project ICT-31. Furthermore, we would like to thank Diego Davila for providing the STAR-CCM+ model of the reactor and corresponding simulation data as well as Peter Mas and Daniel Berger for the valuable feedback and discussions.

\bibliographystyle{plain}
\bibliography{references}


\newpage
\appendix

\section{Reduced Order Model Matrices}\label{app:ROMmatrices}

  \tiny
  \begin{equation*}\label{eq:explicitROM_PNL}
    \bm{P}_1= 10 \left(\begin{array}{rrrrrrrr}
      -0.0840 & -0.0617 & -0.0127 &  0.0671 & -0.2850 &  0.5396 & -0.1903 & -0.3869\\
      -0.0890 &  0.6296 &  2.9940 &  6.6855 &  7.7643 &  5.1206 & -4.3235 &  7.9500\\
      -0.0861 &  0.2482 &  0.2004 & -0.2318 & -0.1732 &  0.0256 & -0.4423 &  0.1641\\
      -0.0851 &  0.1043 & -0.2083 &  0.0653 &  0.2102 & -0.0036 &  0.1797 & -0.3952\\
      -0.0878 &  0.4826 &  1.7418 &  2.7261 &  1.5563 &  0.4511 & -0.5013 & -2.0087\\
      -0.0844 &  0.0056 & -0.1620 &  0.2623 & -0.2392 & -0.1433 & -0.0309 &  0.3140\\
      -0.0865 &  0.3025 &  0.7559 &  1.2072 &  0.9995 &  0.6310 &  0.3663 &  0.7116\\
      -0.0848 &  0.0570 & -0.0551 &  0.1500 & -0.1855 & -0.0130 & -0.0766 &  0.0654\\
    \end{array}\right)
  \end{equation*}
  \normalsize

  \tiny
  \begin{equation*}\label{eq:explicitROM_PPsi}
    \bm{P}_2 = 0.1 \left(\begin{array}{rrrrrrrr}
      -0.0588 & -0.0049 & -0.8234 & -1.2021 &  0.0206 & -0.0920 & -0.0310 & -1.0607\\
      -0.1916 & -0.0053 &  1.7155 &  1.6733 &  0.0195 & -0.5895 &  0.0978 &  1.8271\\
      -0.2016 & -0.0457 &  0.5092 & -1.6390 &  0.2257 & -0.7922 &  0.1433 & -0.0804\\
       0.1373 &  0.0001 & -0.4994 & -0.1023 &  0.1582 &  1.3724 &  0.0653 & -0.0457\\
      -0.3212 &  0.0522 & -0.3326 &  0.8046 & -0.0687 & -0.8529 &  0.0308 & -0.3995\\
       0.5223 &  0.0403 &  0.2901 &  0.4714 & -0.0970 & -0.4043 &  0.0496 &  0.0920\\
      -0.4669 & -0.0274 &  0.0518 &  0.3782 &  0.0708 & -0.1312 &  0.0059 &  0.0796\\
      -0.1423 &  0.0351 &  0.2112 & -0.3527 & -0.0910 &  0.3584 &  0.0065 & -0.1120\\  
    \end{array}\right)
  \end{equation*}
  \normalsize

  \tiny
  \begin{equation*}\label{eq:explicitROM_A}
    \bm{A}= 0.001 \left(\begin{array}{rrrrrrrr}
      -0.0033 &  0.5106 & -0.0805 & -0.3398 &  0.6184 & -0.6086 & -0.1324 & -0.0030\\
       0.0044 & -0.6860 &  0.3958 &  0.3218 & -0.1290 &  0.3695 & -0.1028 &  0.0389\\
      -0.0020 &  0.3205 & -1.3243 &  0.1441 &  0.2398 & -0.0341 &  0.1491 &  0.1299\\
       0.0011 & -0.1810 &  0.5726 & -0.4061 &  0.0250 &  0.0398 &  0.0258 & -0.0663\\
      -0.0011 &  0.1767 & -0.4966 &  0.4591 & -0.0735 &  0.0108 & -0.0240 & -0.0683\\
       0.0010 & -0.1496 &  0.4312 & -0.2221 &  0.1372 & -0.0639 & -0.0404 &  0.0890\\
      -0.0001 &  0.0157 & -0.0392 &  0.0550 &  0.0039 & -0.0619 & -0.1744 &  0.2296\\
      -0.0001 &  0.0144 & -0.0460 &  0.0337 &  0.0115 & -0.0088 & -0.0152 & -0.0318\\
    \end{array}\right)
  \end{equation*}
  \normalsize

\end{document}